\documentclass[11pt]{article}
\usepackage{amscd}
\usepackage{amsfonts}
\usepackage{amsmath}
\usepackage{amssymb}
\usepackage{amsthm}
\usepackage{bbm}
\usepackage{CJK}
\usepackage{fancyhdr}
\usepackage{graphicx}
\usepackage{hyperref}
\usepackage{indentfirst}
\usepackage{latexsym}
\usepackage{mathrsfs}
\usepackage{xypic}
\usepackage{amsmath,amssymb,amsfonts,amsthm,fancyhdr}
\usepackage{epsfig,graphicx,picins,picinpar,subfigure}
\usepackage{pstricks}
\usepackage{amssymb}
\usepackage{xypic}
\usepackage[compress]{cite}

\newtheorem{theorem}{Theorem}[section]
\newtheorem{prop}[theorem]{Proposition}
\newtheorem{lemma}[theorem]{Lemma}
\newtheorem{coro}[theorem]{Corollary}
\newtheorem{prop-def}{Proposition-Definition}[section]

\theoremstyle{definition}
\newtheorem{defn}[theorem]{Definition}

\newtheorem{remark}[theorem]{Remark}
\newtheorem{exam}[theorem]{Example}

\usepackage[top=1in,bottom=1in,left=1.25in,right=1.25in]{geometry}
\def\<{\langle}
\def\>{\rangle}

\date{\today}
\begin{document}
\renewcommand{\baselinestretch}{1.2}
\renewcommand{\arraystretch}{1.0}
\title{\bf  Nonabelian embedding tensors on 3-Lie algebras and 3-Leibniz-Lie algebras}
\author{{\bf Wen Teng$^{1}$,   Xiansheng Dai$^{2}$}\\
{\small 1. School of Mathematics and Statistics, Guizhou University of Finance and Economics} \\
{\small  Guiyang  550025, P. R. of China}\\
  {\small E-mail: tengwen@mail.gufe.edu.cn (Wen Teng)} \\
{\small 2.   School of Mathematical Sciences, Guizhou Normal University}\\
{\small Guizhou  $550005$, P. R. of China}\\
 {\small E-mail:daisheng158@126.com (Xiansheng Dai)}}
 \maketitle
\begin{center}
\begin{minipage}{13.cm}

{\bf Abstract} In this paper, first we introduce the notion of a nonabelian embedding tensor on the 3-Lie algebra. Then, we introduce the notion
of a 3-Leibniz-Lie algebra, which is the underlying algebraic structure of a nonabelian
embedding tensor on the 3-Lie algebra, and can also be viewed as a nonabelian generalization of a 3-Leibniz algebra. Next we develop   the cohomology of    nonabelian embedding tensors on 3-Lie algebras  with coefficients in a suitable representation  and use the first cohomology group to characterize infinitesimal deformations. Finally, we investigate nonabelian embedding tensors  on 3-Lie algebras induced by  Lie algebras.

 \smallskip

 \smallskip

 \end{minipage}
 \end{center}
 \normalsize\vskip0.5cm

\section{Introduction}
\def\theequation{\arabic{section}. \arabic{equation}}
\setcounter{equation} {0}

The concept of  embedded tensor \cite{Nicolai} and related tensor hierarchies provide a useful tool for constructing the theory of supergravity and higher gauge theory\cite{Bergshoeff}.
See \cite{Bonezzi1,Bonezzi2,deWit1,deWit2,deWit3,deWit4,Das,Das1} and the references therein for a great deal of literature on embedding tensors and  related tensor hierarchies.
In \cite{ Kotov}, the authors first observed the mathematical essence behind the embedding tensor and proved that the embedding tensor naturally produced Leibniz algebra.
In the application of physics, they observed that in the construction of the corresponding gauge theory, they focused more on Leibniz algebra than embedding tensor.
In \cite{Sheng}, Sheng, Tang and Zhu considered cohomology, deformations and homotopy theory for embedding tensors and Lie-Leibniz triples.
Later on,  the deformation and cohomology theory of embedding tensors on 3-Lie algebras are given in  \cite{Hu}.   Tang and Sheng \cite{Tang} first proposed the nonabelian embedding tensor on Lie algebras, which is a nonabelian generalization of the embedding tensor, and gave the algebraic structure behind the nonabelian embedding tensor as Leibniz-Lie algebras. Furthermore, the nonabelian embedding tensor on Lie algebras has been extended to the Hom setting in \cite{Teng}.

On the other hand, Filippov \cite{Filippov} first introduced the concepts of 3-Lie algebras and more generally, $n$-Lie algebras (also called Filippov algebras).
In recent years, 3-Lie algebra has been widely studied and applied in the fields of mathematics and physics, including string theory, Nambu mechanics and  M2-branes \cite{Nambu,Bagger,Ho,Gustavsson}.
Further research on 3-Lie algebras could be found in \cite{Kasymov,Takhtajan, Arfa,Mignel, Liu, Liu1,LiuW,Sheng1, Xu, Zhang,Zhao,Hou,Bai,Teng1} and references cited therein.
Motivated by Tang's  \cite{Tang} terminology of nonabelian embedding tensors and considering the importance of 3-Lie algebras,  cohomology and deformation theories, we mainly study the nonabelian embedding tensors on 3-Lie algebras in this paper.

  This paper is organized as follows.
  Section 2  first recalls some basic notions of 3-Lie algebras and 3-Leibniz algebras.   Then we introduce the coherent  action  of a 3-Lie algebra  on another 3-Lie algebra  and   the notion of nonabelian embedding tensors on  3-Lie algebras with respect to a coherent action .
In Section 3, the notion of 3-Leibniz-Lie algebra is introduced as the basic algebraic structure of a nonabelian embedding tensor on the 3-Lie algebra. Naturally, a 3-Leibniz-Lie algebra induces a 3-Leibniz algebra.
In Section 4, the cohomology theory of nonabelian embedding tensors on 3-Lie algebras is introduced.  As   application, we characterize the infinitesimal deformation   using the first cohomology group.
In Section 5, we investigate nonabelian embedding tensors on 3-Lie algebras induced by Lie algebras.

 All vector spaces and algebras considered in this paper are on the field $\mathbb{K}$ with the characteristic of 0.

\section{Nonabelian embedding tensors on 3-Lie algebras}
\def\theequation{\arabic{section}.\arabic{equation}}
\setcounter{equation} {0}

This section  recalls some basic notions of 3-Lie algebras and 3-Leibniz algebras. After that, we introduce the coherent  action  of a 3-Lie algebra  on another 3-Lie algebra, and we introduce the concept of nonabelian embedding tensors on  3-Lie algebras  by  its coherent  action  as a nonabelian generalization of embedding tensors on 3-Lie algebras \cite{Hu}.

\begin{defn}(see\cite{Filippov})  A  3-Lie algebra    is a  pair $(L, [-, -, -]_L)$ consisting of a vector space  $L$  and   a   skew-symmetric trilinear map
 $[-, -, -]_L: \wedge^3 L\rightarrow L$ such that
\begin{align}
&[l_1,l_2,[l_3,l_4,l_5]_L]_L=[[l_1,l_2,l_3]_L,l_4,l_5]_L+[l_3,[l_1,l_2,l_4]_L,l_5]_L+[l_3,l_4,[l_1,l_2,l_5]_L]_L,  \label{2.1}
\end{align}
for any  $l_i\in L$.
\end{defn}

A   homomorphism   between two 3-Lie algebras $(L_1, [-,-,-]_{L_1})$ and $(L_2, [-,-,-]_{L_2})$ is a linear map $f: L_1\rightarrow L_2$ satisfies
$f([l_1,l_2,l_3]_{L_1})=[f(l_1),f(l_2),f(l_3)]_{L_2}, \ \forall l_1,l_2,l_3\in L_1.$

\begin{defn}
(1) (see \cite{Kasymov}) A representation of a 3-Lie algebra  $(L, [ -, -, - ]_L)$  on a vector space  $H$  is a   linear map $\rho:  \wedge^2 L \rightarrow \mathrm{End}(H)$, such that
\begin{align}
& \rho([l_1, l_2,  l_3]_L,l_4)=\rho(l_2,  l_3)\rho(l_1,l_4)+\rho(l_3,  l_1)\rho(l_2,  l_4)+\rho(l_1,  l_2)\rho(l_3,  l_4),\label{2.2}\\
& \rho(l_1,  l_2)\rho(l_3,  l_4)=\rho(l_3,  l_4)\rho(l_1,  l_2)+ \rho([l_1, l_2,  l_3]_L,l_4)+ \rho(l_3,[l_1, l_2,  l_4]_L),\label{2.3}\
\end{align}
for all  $l_1, l_2,  l_3,l_4\in L$.  We also denote a representation of $ L $  on   $H$  by  $(H;  \rho)$. \\
(2)  A coherent  action of  a  3-Lie algebra  $(L, [ -, -, - ]_L)$  on another 3-Lie algebra  $(H, [ -,-,- ]_H)$  is a   linear map $\rho:  \wedge^2 L \rightarrow \mathrm{End}(H)$ satisfying Eqs.  \eqref{2.2}, \eqref{2.3} and
\begin{align}
\rho(l_1,l_2)[h_1,h_2,h_3]_H=&[\rho(l_1,l_2)h_1,h_2,h_3]_H+[h_1,\rho(l_1,l_2)h_2,h_3]_H+[h_1,h_2,\rho(l_1,l_2)h_3]_H,\label{2.4}\\
[\rho(l_1,l_2)h_1,h_2,h_3]_H=&0, \label{2.5}
\end{align}
for all $l_1,l_2,l_3\in L$ and $h_1,h_2,h_3\in H$. We denote a coherent  action of  $ L $  on   $H$  by $(H,[-,-,-]_H; \rho^{\dag}).$
\end{defn}

\begin{exam}
 Let $(H, [-,-,- ]_H)$ be a  3-Lie algebra. Define $ad:  \wedge^2 H \rightarrow \mathrm{End}(H)$ by
$ad(h_1,h_2)(h):=[h_1,h_2,h], \forall \ h_1,h_2,h\in H.$
Then $(H;ad)$ is a representation of $(H, [-,-,- ]_H)$, which is called the adjoint representation. Furthermore, if $ad$ satisfies
$$[ad(h_1,h_2)h'_1,h'_2,h'_3]_H=0,  \forall \ h'_1,h'_2,h'_3\in H,$$
then $(H, [-,-,-]_H; ad^{\dag})$ is a coherent adjoint  action  of $(H, [-,-,-]_H)$.
\end{exam}

\begin{defn} (see\cite{Casas})
A   3-Leibniz algebra   is a vector space $\mathcal{L}$ together with a
 trilinear operation $[-,-,- ]_{\mathcal{L}}: \mathcal{L}\otimes \mathcal{L}\otimes \mathcal{L}  \rightarrow \mathcal{L}$   such that
\begin{align*}
&[l_1,l_2,[l_3,l_4,l_5]_\mathcal{L}]_\mathcal{L}=[[l_1,l_2,l_3]_\mathcal{L},l_4,l_5]_\mathcal{L}+[l_3,[l_1,l_2,l_4]_\mathcal{L},l_5]_\mathcal{L}+[l_3,l_4,[l_1,l_2,l_5]_\mathcal{L}]_\mathcal{L},
\end{align*}
for any  $l_i\in \mathcal{L}$.
\end{defn}

\begin{prop}
Let $(L, [-,-,- ]_L)$ and $(H, [-,-,- ]_H)$ be two 3-Lie algebras  and  $\rho:  \wedge^2 L \rightarrow \mathrm{End}(H)$ a bilinear map.   Then $(H, [-,-,-]_H; \rho^{\dag})$ is  a  coherent   action  of  $L$ if and only if
 $L \oplus H$ is a 3-Leibniz algebra under the following map:
\begin{eqnarray*}
&[l_1+h_1, l_2+h_2,l_3+h_3]_{\rho}:=[l_1, l_2, l_3]_L+\rho(l_1,l_2)h_3+[h_1,h_2,h_3]_H,
\end{eqnarray*}
for any    $l_1, l_2, l_3\in L$ and $h_1,h_2,h_3\in H$.  $(L \oplus H, [-,-,-]_{\rho})$  is called
the   nonabelian hemisemidirect product 3-Leibniz algebra, and denoted by $L\ltimes_{\rho}H$.
\end{prop}

\begin{proof}  For all  $l_1, l_2, l_3, l_4, l_5 \in L$ and $h_1,h_2,h_3,h_4,h_5\in H$,  by Eqs. \eqref{2.1},  \eqref{2.3}, \eqref{2.4} and \eqref{2.5}, we have
\begin{align*}
&[l_1+h_1, l_2+h_2,[l_3+h_3,l_4+h_4, l_5+h_5]_{\rho}]_{\rho}-[[l_1+h_1, l_2+h_2,l_3+h_3]_{\rho},l_4+h_4, l_5+h_5]_{\rho}\\
&-[l_3+h_3,[l_1+h_1, l_2+h_2,l_4+h_4]_{\rho}, l_5+h_5]_{\rho}-[l_3+h_3,l_4+h_4,[l_1+h_1, l_2+h_2, l_5+h_5]_{\rho}]_{\rho}\\
=&[l_1, l_2,[l_3, l_4, l_5]_L]_L+\rho(l_1,l_2)\rho(l_3,l_4)h_5+\rho(l_1,l_2)[h_3,h_4,h_5]_H+[h_1,h_2,\rho(l_3,l_4)h_5]_H\\
&+[h_1,h_2,[h_3,h_4,h_5]_H]_H-[[l_1, l_2, l_3]_L,l_4,l_5]_L-\rho([l_1, l_2, l_3]_L,l_4)h_5-[\rho(l_1,l_2)h_3, h_4,  h_5]_H\\
&-[[h_1,h_2,h_3]_H, h_4,  h_5]_H-[l_3,[l_1, l_2, l_4]_L,l_5]_L-\rho(l_3,[l_1, l_2, l_4]_L)h_5-[h_3,\rho(l_1,l_2)h_4,h_5]_H\\
&-[h_3,[h_1,h_2,h_4]_H,  h_5]_{H}-[l_3,l_4,[l_1, l_2, l_5]_L]_L-\rho(l_3,l_4)\rho(l_1,l_2)h_5-\rho(l_3,l_4)[h_1,h_2,h_5]_H\\
&-[h_3,h_4, \rho(l_1,l_2)h_5]_{H}-[h_3,h_4,[h_1,h_2,h_5]_H]_{H}\\
=&[h_1,h_2,\rho(l_3,l_4)h_5]_H-\rho(l_3,l_4)[h_1,h_2,h_5]_H\\
=&0.
\end{align*}
Thus,   $(L \oplus H, [-,-,- ]_{\rho})$   is a 3-Leibniz algebra.

The converse can be proved similarly. We omit the details.
\end{proof}

\begin{defn} (1) A   nonabelian embedding tensor   on  a 3-algebra  $(L, [-,-,- ]_L)$   with respect to  a coherent  action  $(H, [-, -, - ]_H;  \rho^{\dag})$ is a linear map $\Lambda: H\rightarrow L$
satisfying the following equation:
\begin{align}
[\Lambda h_1,\Lambda h_2, \Lambda h_3]_L=&\Lambda(\rho(\Lambda h_1,\Lambda h_2)h_3+[h_1,h_2,h_3]_H), \label{2.6}
\end{align}
for any  $h_1,h_2,h_3\in H$.

(2) A    nonabelian embedding tensor 3-Lie algebra  is a triple $(H,L,\Lambda)$ consisting of a   3-Lie algebra   $(L, [-,-,-]_L)$,  a coherent  action  $(H,[-,-,-]_H; \rho^{\dag})$ of $L$  and a nonabelian embedding tensor $\Lambda: H\rightarrow L$.  We denote a nonabelian embedding tensor 3-Lie algebra  $(H,L,\Lambda)$ by the notation $H\stackrel{\Lambda}{\longrightarrow}L$.

(3) Let $H\stackrel{\Lambda_1}{\longrightarrow}L$ and $H\stackrel{\Lambda_2}{\longrightarrow}L$  be two nonabelian embedding tensor 3-Lie algebras.  Then a   homomorphism  from $H\stackrel{\Lambda_1}{\longrightarrow}L$ to $H\stackrel{\Lambda_2}{\longrightarrow}L$  consists of two   3-Lie algebras
homomorphisms $f_L:L\rightarrow L$ and  $f_{H}:H\rightarrow H$ satisfying the following equations
\begin{align}
\Lambda_2\circ f_H=&f_L\circ \Lambda_1,\label{2.7}\\
f_{H}(\rho(l_1,l_2)h)= &\rho(f_L(l_1),f_L(l_2))f_{H}(h).\label{2.8}
\end{align}
for $l_1,l_2\in L$ and $h\in H.$  Furthermore, if   $f_L$ and $f_{H}$ are nondegenerate, $(f_L,f_{H})$ is called an isomorphism from $H\stackrel{\Lambda_1}{\longrightarrow}L$ to $H\stackrel{\Lambda_2}{\longrightarrow}L$.
\end{defn}

\begin{remark}
If  $(H, [-,-,- ]_H)$ is an abelian 3-Lie algebra, then we can get that $\Lambda$ is an  embedding tensor on 3-Lie algebra (see \cite{Hu}). In addition, If $\rho=0$, then $\Lambda$ is a 3-Lie algebra homomorphism from $H$ to $L$.
\end{remark}

\begin{exam}  \label{exe:2.8}
Let $H$ be a 4-dimensional linear space spanned by  $\alpha_1,\alpha_2,\alpha_3$ and $\alpha_4$.  We define a
skew-symmetric trilinear map  $[-,-,-]_H:\wedge^3 H\rightarrow H$  by
$$[\alpha_1,\alpha_2,\alpha_3]_H=\alpha_4.$$
Then   $(H, [-,-,-]_H)$  is a  3-Lie algebra. It is obvious that
 $(H, [-,-,-]_H; ad^{\dag})$ is a coherent adjoint  action  of $(H, [-,-,-]_H)$.
  Moreover, for $k\in \mathbb{K}$,
  $$\Lambda=\left(
                      \begin{array}{cccc}
                       1 & 0 & 0 & 0\\
                       0 & 1 & 0 & 0\\
                      0 & 0 & 2k & 0\\
                      0 & 0 & 0 & k\\
                      \end{array}
                    \right)$$
  is a nonabelian embedding tensor   on  the 3-algebra  $(H, [-,-,-]_H)$   with respect to  the  coherent adjoint  action $(H, [-,-,-]_H; ad^{\dag})$.

\end{exam}

Next we use graphs to describe nonabelian embedding tensors on 3-Lie algebras.
\begin{theorem}
 A linear map $\Lambda: H\rightarrow L$ is a nonabelian embedding tensor  on  a 3-Lie algebra    $(L, [-,-,-]_L)$  with respect to   the   coherent  action  $(H,[-,-,-]_H; \rho^{\dag})$  if and only if the graph $Gr(\Lambda)=\{\Lambda h+h~|~h\in H\}$ is a
subalgebra of the nonabelian hemisemidirect product 3-Leibniz algebra $L\ltimes_{\rho}H$.
\end{theorem}

\begin{proof}  Let $\Lambda: H\rightarrow L$ be a linear map. Then for all  $h_1,h_2,h_3\in H$,  we have
\begin{align*}
&[\Lambda h_1+h_1,\Lambda h_2+h_2, \Lambda h_3+h_3]_{\rho}=[\Lambda h_1, \Lambda h_2, \Lambda h_3]_L+\rho(\Lambda h_1, \Lambda h_2)h_3+[h_1,h_2,h_3]_H,
\end{align*}
Thus,    the graph $Gr(\Lambda)=\{\Lambda h+h ~|~ h\in H\}$ is a subalgebra of the  nonabelian  hemisemidirect product
3-Leibniz algebra $L\ltimes_{\rho}H$ if and only if $\Lambda$  satisfies  Eq. \eqref{2.6}, which implies
that $\Lambda$ is a  nonabelian embedding tensor  on     $L$  with respect to   the  coherent  action  $(H,[-,-,-]_H; \rho^{\dag})$.
\end{proof}

Because $H$ and $Gr(\Lambda)$ are isomorphic as linear spaces,  there is an induced 3-Leibniz algebra
structure on $H$.

\begin{coro} \label{cor:3-Leibniz}
Let $H\stackrel{\Lambda}{\longrightarrow}L$    be a nonabelian embedding tensor 3-Lie algebra. If a linear map $[-,-,-]_\Lambda: \wedge^3 H\rightarrow H$ is given by
\begin{align}
&[h_1,h_2,h_3]_\Lambda=\rho(\Lambda h_1, \Lambda h_2)h_3+[h_1,h_2,h_3]_H, \label{2.9}
\end{align}
for any     $h_1,h_2,h_3\in H$. Then $(H,[-,-,-]_\Lambda)$ is a 3-Leibniz algebra. Moreover, $\Lambda$ is a homomorphism from the 3-Leibniz algebra
  $(H,[-,-,-]_\Lambda)$ to the 3-Lie algebra $(L,[-,-,-]_L)$. The 3-Leibniz algebra $(H,[-,-,-]_\Lambda)$ is called the   descendent 3-Leibniz algebra.
\end{coro}

\begin{prop}\label{pro:2.9}
Let   $(f_L,f_{H})$  be a homomorphism from $H\stackrel{\Lambda_1}{\longrightarrow}L$ to $H\stackrel{\Lambda_2}{\longrightarrow}L$. Then $f_{H}$  is a
homomorphism of descendent 3-Leibniz algebra  from $(H,[-,-,-]_{\Lambda_1})$ to $(H,[-,-,-]_{\Lambda_2})$.
\end{prop}

\begin{proof}  For all $h_1,h_2,h_3\in H,$ by Eqs.   \eqref{2.7}, \eqref{2.8} and \eqref{2.9}, we have
\begin{align*}
f_{H}([h_1,h_2,h_3]_{\Lambda_1})=&f_{H}(\rho(\Lambda_1 h_1, \Lambda_1 h_2)h_3+[h_1,h_2,h_3]_H)\\
=&\rho(f_L(\Lambda_1 h_1), f_L(\Lambda_1 h_2))f_{H}(h_3)+f_{H}([h_1,h_2,h_3]_H)\\
=&\rho(\Lambda_2 f_L(h_1), \Lambda_2 f_L(h_2))f_{H}(h_3)+[f_{H}(h_1),f_{H}(h_2),f_{H}(h_3)]_H\\
=&[f_{H}(h_1),f_{H}(h_2),f_{H}(h_3)]_{\Lambda_2}.
\end{align*}
 The proof is finished.
 \end{proof}

\section{3-Leibniz-Lie algebras}
\def\theequation{\arabic{section}.\arabic{equation}}
\setcounter{equation} {0}

In this section, we introduce the concept of 3-Leibniz-Lie algebra as the basic algebraic structure of nonabelian embedding tensor 3-Lie algebra.

\begin{defn}
A 3-Leibniz-Lie algebra  $(H, [-,-,-]_H, \{-,-,-\}_H)$ consists of a 3-Lie algebra  $(H, [-,-,-]_H)$ and a
trilinear product $\{-,-,-\}_H: \wedge^3 H\rightarrow H$   such that
\begin{align}
&\{h_1,h_2,\{h_3,h_4,h_5\}_H\}_H=\{\{h_1,h_2,h_3\}_H,h_4,h_5\}_H+\{h_3,\{h_1,h_2,h_4\}_H, h_5\}_H\nonumber\\
&\ \ \ \ +\{h_3,h_4,\{h_1,h_2,h_5\}_H\}_H+\{[h_1,h_2,h_3]_H,h_4,h_5\}_H+\{h_3,[h_1,h_2,h_4]_H, h_5\}_H,\label{3.1}\\
&\{h_1,h_2,[h_3,h_4,h_5]_H\}_H=[\{h_1,h_2,h_3\}_H,h_4,h_5]_H=0, \label{3.2}
\end{align}
for any $h_1,h_2,h_3,h_4,h_5\in H$.
\end{defn}

A   homomorphism  between two 3-Leibniz-Lie algebras $(H_1, [-,-,-]_{H_1}, \{-,-,-\}_{H_1})$  and  $(H_2, [-,-,-]_{H_2}, \{-,-,-\}_{H_2})$  is a 3-Lie algebra  $f: (H_1,, [-,-,-]_{H_1})\rightarrow (H_2,, [-,-,-]_{H_2}) $ such that
$f(\{h_1,h_2,h_3\}_{H_1})=\{f(h_1),f(h_2),f(h_3)\}_{H_2},\ \ \forall  h_1,h_2,h_3\in H_1.$

\begin{remark}
A 3-Leibniz  algebra $(H,   \{-,-,-\}_H)$ is naturally a   3-Leibniz-Lie algebra  if the 3-Lie algebra  $(H, [-,-,-]_H)$  is abelian.
\end{remark}

\begin{exam} Let $(H, [-,-,-]_H)$ be a 4-dimensional  3-Lie algebra  given in  Example  \ref{exe:2.8}.
We define a trilinear product $\{-,-,-\}_H: \wedge^3 H\rightarrow H$  by
\begin{align*}
&\{\alpha_{i_1},\alpha_{i_2},\alpha_{i_3}\}_H=(-1)^{i_1+i_2+i_3}\alpha_4, \ \ \  i_1,i_2,i_3\in \{1,2,3\},\\
&\{\alpha_4,\alpha_{j_1},\alpha_{j_2}\}_H=\{\alpha_{j_1},\alpha_4,\alpha_{j_2}\}_H=\{\alpha_{j_1},\alpha_{j_2},\alpha_4\}_H=0,\ j_1,j_2\in\{1,2,3,4\}.
\end{align*}
Then  $(H, [-,-,-]_H, \{-,-,-\}_H)$  is a   3-Leibniz-Lie algebra.
\end{exam}

 The following theorem shows that a 3-Leibniz-Lie algebra
naturally induces a 3-Leibniz algebra.

\begin{theorem}
  Let $(H, [-,-,-]_H, \{-,-,-\}_H)$ be a 3-Leibniz-Lie algebra. Then the trilinear product $\langle-,-,-\rangle_H: \wedge^3 H\rightarrow H$ given by
\begin{align}
\langle h_1,h_2,h_3\rangle_H:=[h_1,h_2,h_3]_H+\{h_1,h_2,h_3\}_H, \label{3.3}
\end{align}
for any $h_1,h_2,h_3\in H,$ defines a 3-Leibniz  algebra structure on $H$, which is denoted by $(H,\langle-,-,-\rangle_H)$
and called the  subadjacent 3-Leibniz  algebra.
\end{theorem}

\begin{proof} For any $h_1,h_2,h_3,h_4,h_5\in H$, by Eqs. \eqref{2.1}, \eqref{3.1}, \eqref{3.2} and \eqref{3.3}, we have
\begin{align*}
&\langle h_1,h_2,\langle h_3,h_4,h_5 \rangle_H\rangle_H-\langle \langle h_1,h_2, h_3\rangle_H,h_4,h_5\rangle_H-\langle h_3,\langle h_1,h_2,h_4\rangle_H,h_5\rangle_H\\
&-\langle h_3,h_4,\langle h_1,h_2,h_5 \rangle_H\rangle_H\\
=&[ h_1,h_2,[h_3,h_4,h_5]_H]_H+[ h_1,h_2,\{h_3,h_4,h_5\}_H]_H+\{h_1,h_2,[h_3,h_4,h_5]_H\}_H\\
&+\{h_1,h_2,\{h_3,h_4,h_5\}_H\}_H-[ [h_1,h_2,h_3]_H,h_4,h_5]_H-[\{h_1,h_2,h_3\}_H,h_4,h_5]_H\\
&-\{[h_1,h_2,h_3]_H,h_4,h_5\}_H-\{\{h_1,h_2,h_3\}_H,h_4,h_5\}_H-[ h_3,[h_1,h_2,h_4]_H,h_5]_H\\
&-[h_3,\{h_1,h_2,h_4\}_H,h_5]_H-\{h_3,[h_1,h_2,h_4]_H,h_5\}_H-\{h_3,\{h_1,h_2,h_4\}_H,h_5\}_H\\
&-[ h_3,h_4,[h_1,h_2,h_5]_H]_H-[ h_3,h_4,\{h_1,h_2,h_5\}_H]_H-\{h_3,h_4,[h_1,h_2,h_5]_H\}_H\\
&-\{ h_3,h_4,\{h_1,h_2,h_5\}_H\}_H\\
=&\{h_1,h_2,\{h_3,h_4,h_5\}_H\}_H-\{[h_1,h_2,h_3]_H,h_4,h_5\}_H-\{\{h_1,h_2,h_3\}_H,h_4,h_5\}_H\\
&-\{h_3,[h_1,h_2,h_4]_H,h_5\}_H-\{h_3,\{h_1,h_2,h_4\}_H,h_5\}_H-\{ h_3,h_4,\{h_1,h_2,h_5\}_H\}_H\\
=&0.
\end{align*}
Hence, $(H,\langle-,-,-\rangle_H)$  is a 3-Leibniz  algebra.
\end{proof}

 The following theorem shows that a nonabelian embedding tensor 3-Lie algebra
 induces a 3-Leibniz-Lie algebra.

\begin{theorem}
 Let $H\stackrel{\Lambda}{\longrightarrow}L$    be a nonabelian embedding tensor 3-Lie algebra.  Then  $(H, [-,-,-]_H, \{-,-,-\}_{\Lambda})$ is a 3-Leibniz-Lie algebra, where
\begin{align}
\{h_1,h_2,h_3\}_\Lambda:=\rho(\Lambda h_1,\Lambda h_2)h_3, \label{3.4}
\end{align}
for any $h_1,h_2,h_3\in H$.
\end{theorem}

\begin{proof} For all $h_1,h_2,h_3,h_4,h_5\in H$, by Eqs. \eqref{2.3}, \eqref{2.6}  and \eqref{3.4}, we have
\begin{align*}
&\{\{h_1,h_2,h_3\}_\Lambda,h_4,h_5\}_\Lambda+\{h_3,\{h_1,h_2,h_4\}_\Lambda, h_5\}_\Lambda +\{h_3,h_4,\{h_1,h_2,h_5\}_\Lambda\}_\Lambda\\
&+\{[h_1,h_2,h_3]_H,h_4,h_5\}_\Lambda+\{h_3,[h_1,h_2,h_4]_H,h_5\}_\Lambda-\{h_1,h_2,\{h_3,h_4,h_5\}_\Lambda\}_\Lambda\\
=&\rho(\Lambda\rho(\Lambda h_1,\Lambda h_2)h_3,\Lambda h_4)h_5+\rho(\Lambda h_3,\Lambda\rho(\Lambda h_1,\Lambda h_2)h_4) h_5+\rho(\Lambda h_3,\Lambda h_4)\rho(\Lambda h_1,\Lambda h_2)h_5\\
&+\rho(\Lambda[h_1,h_2,h_3]_H,\Lambda h_4)h_5+\rho(\Lambda h_3,\Lambda[h_1,h_2,h_4]_H) h_5-\rho(\Lambda h_1,\Lambda h_2)\rho(\Lambda h_3,\Lambda h_4)h_5
\end{align*}
\begin{align*}
=&\rho(\Lambda\rho(\Lambda h_1,\Lambda h_2)h_3,\Lambda h_4)h_5+\rho(\Lambda h_3,\Lambda\rho(\Lambda h_1,\Lambda h_2)h_4) h_5+\rho(\Lambda h_3,\Lambda h_4)\rho(\Lambda h_1,\Lambda h_2)h_5\\
&+\rho([\Lambda h_1,\Lambda h_2,\Lambda h_3]_L-\Lambda\rho(\Lambda h_1,\Lambda h_2)h_3,\Lambda h_4)h_5+\rho(\Lambda h_3,[\Lambda h_1,\Lambda h_2,\Lambda h_4]_L\\
&-\Lambda\rho(\Lambda h_1,\Lambda h_2)h_4) h_5-\rho(\Lambda h_1,\Lambda h_2)\rho(\Lambda h_3,\Lambda h_4)h_5\\
=&\rho(\Lambda h_3,\Lambda h_4)\rho(\Lambda h_1,\Lambda h_2)h_5+\rho([\Lambda h_1,\Lambda h_2,\Lambda h_3]_L,\Lambda h_4)h_5+\rho(\Lambda h_3,[\Lambda h_1,\Lambda h_2,\Lambda h_4]_L) h_5\\
&-\rho(\Lambda h_1,\Lambda h_2)\rho(\Lambda h_3,\Lambda h_4)h_5\\
=&0.
\end{align*}
Furthermore, by Eqs. \eqref{2.4} and  \eqref{2.5}, we have
\begin{align*}
&\{h_1,h_2,[h_3,h_4,h_5]_H\}_\Lambda=\rho(\Lambda h_1,\Lambda h_2)[h_3,h_4,h_5]_H=0,\\
&[\{h_1,h_2,h_3\}_\Lambda,h_4,h_5]_H=[\rho(\Lambda h_1,\Lambda h_2)h_3,h_4,h_5]_H=0.
\end{align*}
Thus,    $(H, [-,-,-]_H, \{-,-,-\}_\Lambda)$ is a 3-Leibniz-Lie algebra.
\end{proof}

\begin{prop}
Let   $(f_L,f_{H})$  be a homomorphism from $H\stackrel{\Lambda_1}{\longrightarrow}L$ to $H\stackrel{\Lambda_2}{\longrightarrow}L$.  Then $f_{H}$  is a
homomorphism of 3-Leibniz-Lie algebras from $(H, [-,-,-]_H, \{-,-,-\}_{\Lambda_1})$ to $(H, [-,-,-]_H, \{-,-,-\}_{\Lambda_2})$.
\end{prop}

\begin{proof} For any $h_1,h_2,h_3\in H$, by Eqs. \eqref{2.7}, \eqref{2.8} and  \eqref{3.4}, we have
\begin{align*}
f_H(\{h_1,h_2,h_3\}_{\Lambda_1})=&f_H(\rho(\Lambda_1 h_1,\Lambda_1 h_2)h_3)\\
=&\rho(f_L(\Lambda_1 h_1),f_L(\Lambda_1 h_2))f_H(h_3)\\
=&\rho(\Lambda_2 f_H(h_1),\Lambda_2 f_H( h_2))f_H(h_3)\\
=&\{f_H(h_1),f_H(h_2),f_H(h_3)\}_{\Lambda_2}.
\end{align*}
The proof is finished.
\end{proof}

\section{ Cohomology theory  of   nonabelian embedding tensors on 3-Lie algebras}
\def\theequation{\arabic{section}.\arabic{equation}}
\setcounter{equation} {0}

In this section,  we recall some basic results of  representations and cohomologies of
3-Leibniz algebras.  We construct a representation of the descendent 3-Leibniz algebra $(H,[-,-,-]_\Lambda)$
on the vector space $L$,  and define the cohomologies of a  nonabelian embedding tensor on 3-Lie algebras.
As  application, we characterize the infinitesimal deformation   using the first cohomology group.

\begin{defn}  A representation of the 3-Leibniz algebra $(\mathcal{H},[-,-,-]_{\mathcal{H}})$ is a
vector space $V$  equipped with 3 actions
\begin{align*}
\mathfrak{l}:\mathcal{H}\otimes \mathcal{H}\otimes V\rightarrow V,\\
\mathfrak{m}:\mathcal{H}\otimes V\otimes \mathcal{H}\rightarrow V,\\
\mathfrak{r}:V\otimes \mathcal{H}\otimes \mathcal{H}\rightarrow V,
\end{align*}
satisfying for any  $a_1,a_2,a_3,a_4,a_5\in \mathcal{H}$ and $u\in V$
\begin{small}
\begin{align}
&\mathfrak{l}(a_1,a_2,\mathfrak{l}(a_3,a_4,u))=\mathfrak{l}([a_1,a_2,a_3]_{\mathcal{H}}, a_4, u)+ \mathfrak{l}(a_3,[a_1,a_2,a_4]_{\mathcal{H}},u)+\mathfrak{l}(a_3,a_4,\mathfrak{l}(a_1,a_2,u)),\label{4.1}\\
&\mathfrak{l}(a_1,a_2,\mathfrak{m}(a_3,u, a_5))=\mathfrak{m}([a_1,a_2,a_3]_{\mathcal{H}}, u, a_5)+ \mathfrak{m}(a_3,l(a_1,a_2,u),a_5)+\mathfrak{m}(a_3,u,[a_1,a_2,a_5]_{\mathcal{H}}),\label{4.2}\\
&\mathfrak{l}(a_1,a_2,\mathfrak{r}(u,a_4, a_5))=\mathfrak{r}(\mathfrak{l}(a_1,a_2,u), a_4, a_5)+ \mathfrak{r}(u,[a_1,a_2,a_4]_\mathcal{H},a_5)+\mathfrak{r}(u,a_4,[a_1,a_2,a_5]_{\mathcal{H}}),\label{4.3}\\
&\mathfrak{m}(a_1,u,[a_3,a_4, a_5]_{\mathcal{H}})=\mathfrak{r}(\mathfrak{m}(a_1,u,a_3), a_4, a_5)+ \mathfrak{m}(a_3,\mathfrak{m}(a_1,u,a_4),a_5)+\mathfrak{l}(a_3,a_4,\mathfrak{m}(a_1,u,a_5)),\label{4.4}\\
&\mathfrak{r}(u,a_2,[a_3,a_4, a_5]_{\mathcal{H}})=\mathfrak{r}(\mathfrak{r}(u,a_2,a_3), a_4, a_5)+ \mathfrak{m}(a_3,\mathfrak{r}(u,a_2,a_4),a_5)+\mathfrak{l}(a_3,a_4,\mathfrak{r}(u,a_2,a_5)).\label{4.5}
\end{align}
\end{small}
\end{defn}

For $n\geq 1$, denote the $n$-cochains of  3-Leibniz algebra $(\mathcal{H},[-,-,-]_{\mathcal{H}})$  with coefficients in a representation  $(V;\mathfrak{l},\mathfrak{m},\mathfrak{r})$ by
$$\mathcal{C}^n_{\mathrm{3Leib}}(\mathcal{H},V)=\mathrm{Hom}(\overbrace{\wedge^2\mathcal{H}\otimes\cdots \otimes \wedge^2 \mathcal{H}}^{n-1}\otimes \mathcal{H}, V).$$
  The coboundary map $\delta:\mathcal{C}^n_{\mathrm{3Leib}}(\mathcal{H},V)\rightarrow \mathcal{C}^{n+1}_{\mathrm{3Leib}}(\mathcal{H},V)$, for $A_i=a_i\wedge b_i\in \wedge^2 \mathcal{H}, 1\leq i\leq n$ and $c\in \mathcal{H}$, as
\begin{small}
\begin{align*}
&(\delta \varphi)(A_1,A_2, \cdots, A_n,c)\\
=&\sum_{1\leq j<k\leq n}(-1)^j\varphi(A_1,\cdots,\widehat{A_j},\cdots,A_{k-1}, a_k  \wedge[a_j,b_j,b_k]_{\mathcal{H}}+[a_j,b_j,a_k]_{\mathcal{H}}\wedge  b_k,\cdots,A_n,c)\\
&+\sum_{j=1}^n(-1)^j\varphi(A_1,\cdots,\widehat{A_j},\cdots,A_{n},[a_j,b_j,c]_{\mathcal{H}})\\
&+\sum_{j=1}^n(-1)^{j+1}\mathfrak{l}(A_j,\varphi(A_1,\cdots,\widehat{A_j},\cdots,A_{n},c))\\
&+(-1)^{n+1}(\mathfrak{m}(a_n, \varphi(A_1,\cdots,A_{n-1},v_n),c)+\mathfrak{r}(\varphi(A_1,\cdots,A_{n-1},a_n),  b_n, c)).
\end{align*}
\end{small}
It was proved in \cite{Casas, Takhtajan} that $\delta^2=0$. Therefore,
$(\oplus_{n=1}^{+\infty}\mathcal{C}^n_{\mathrm{3Leib}}(\mathcal{H},V),\delta)$ is a cochain complex.

Let $H\stackrel{\Lambda}{\longrightarrow}L$    be a nonabelian embedding tensor 3-Lie algebra.  By Corollary \ref{cor:3-Leibniz}, $(H,[-,-,-]_\Lambda)$ is a  3-Leibniz algebra.  Next we give a representation of $(H,[-,-,-]_\Lambda)$ on $L$.

\begin{lemma}  With above notations.  Define 3 actions
\begin{align*}
\mathfrak{l}_\Lambda:H\otimes H\otimes L\rightarrow L,\mathfrak{m}_\Lambda:H\otimes L\otimes H\rightarrow L,\mathfrak{r}_\Lambda:L\otimes H\otimes H\rightarrow L,
\end{align*}
by
\begin{align*}
\mathfrak{l}_\Lambda(h_1,h_2,l)&=[\Lambda h_1,\Lambda h_2,l]_L,\\
\mathfrak{m}_\Lambda(h_1,l,h_2)&=[\Lambda h_1,l,\Lambda h_2]_L-\Lambda\rho(\Lambda h_1,l)h_2,\\
\mathfrak{r}_\Lambda(l,h_1,h_2)&=[l,\Lambda h_1,\Lambda h_2]_L-\Lambda\rho(l,\Lambda h_1)h_2,
\end{align*}
for any $h_1,h_2\in H,l\in L.$ Then $(L;\mathfrak{l}_\Lambda,\mathfrak{m}_\Lambda,\mathfrak{r}_\Lambda)$ is a representation of the descendent 3-Leibniz algebra $(H,[-,-,-]_\Lambda)$.
\end{lemma}

\begin{proof}  For all $h_1,h_2,h_3,h_4,h_5\in H$ and $l\in L$,
by Eqs. \eqref{2.1}, \eqref{2.3}-\eqref{2.6} and \eqref{2.9}, we have
\begin{align*}
&\mathfrak{l}_\Lambda(h_1,h_2,\mathfrak{l}_\Lambda(h_3,h_4,l))-\mathfrak{l}_\Lambda([h_1,h_2,h_3]_{\Lambda}, h_4, l)- \mathfrak{l}_\Lambda(h_3,[h_1,h_2,h_4]_{\Lambda},l)-\mathfrak{l}_\Lambda(h_3,h_4,\mathfrak{l}_\Lambda(h_1,h_2,l))\\
=&[\Lambda h_1,\Lambda h_2,[\Lambda h_3,\Lambda h_4,l]_L]_L-[[\Lambda h_1, \Lambda h_2, \Lambda h_3]_L, \Lambda h_4, l]_L-[\Lambda h_3,[\Lambda h_1, \Lambda h_2,\Lambda h_4]_L,l]_L\\
&-[\Lambda h_3,\Lambda h_4,[\Lambda h_1,\Lambda h_2,l]_L]_L\\
=&0,\\
&\mathfrak{l}_\Lambda(h_1,h_2,\mathfrak{m}_\Lambda(h_3,l, h_5))-\mathfrak{m}_\Lambda([h_1,h_2,h_3]_{\Lambda}, l, h_5)- \mathfrak{m}_\Lambda(h_3,\mathfrak{l}_\Lambda(h_1,h_2,l),h_5)-\mathfrak{m}_\Lambda(h_3,l,[h_1,h_2,h_5]_{\Lambda})\\
=&[\Lambda h_1,\Lambda h_2,[\Lambda h_3,l,\Lambda h_5]_L]_L-[\Lambda h_1,\Lambda h_2,\Lambda\rho(\Lambda h_3,l)h_5]_L-[[\Lambda h_1,\Lambda h_2,\Lambda h_3]_{L},l,\Lambda h_5]_L\\
&+\Lambda\rho([\Lambda h_1,\Lambda h_2,\Lambda h_3]_{L},l)h_5- [\Lambda h_3,[\Lambda h_1,\Lambda h_2,l]_L,\Lambda h_5]_L+\Lambda\rho(\Lambda h_3,[\Lambda h_1,\Lambda h_2,l]_L)h_5\\
&-[\Lambda h_3,l, [\Lambda h_1,\Lambda h_2,\Lambda h_5]_{L}]_L+\Lambda\rho(\Lambda h_3,l)\rho(\Lambda h_1, \Lambda h_2)h_5+\Lambda\rho(\Lambda h_3,l)[h_1,h_2,h_5]_H\\
=&-[\Lambda h_1,\Lambda h_2,\Lambda\rho(\Lambda h_3,l)h_5]_L+\Lambda\rho([\Lambda h_1,\Lambda h_2,\Lambda h_3]_{L},l)h_5+\Lambda\rho(\Lambda h_3,[\Lambda h_1,\Lambda h_2,l]_L)h_5\\
&+\Lambda\rho(\Lambda h_3,l)\rho(\Lambda h_1, \Lambda h_2)h_5+\Lambda\rho(\Lambda h_3,l)[h_1,h_2,h_5]_H\\
=&- \Lambda \rho(\Lambda h_1,\Lambda h_2)\Lambda\rho(\Lambda h_3,l)h_5-\Lambda[h_1,h_2,\rho(\Lambda h_3,l)h_5]_H+\Lambda\rho(\Lambda h_1, \Lambda h_2)\rho(\Lambda h_3,l)h_5\\
&+\Lambda\rho(\Lambda h_3,l)[h_1,h_2,h_5]_H\\
=&-\Lambda[h_1,h_2,\rho(\Lambda h_3,l)h_5]_H+\Lambda\rho(\Lambda h_3,l)[h_1,h_2,h_5]_H\\
=&\Lambda[\rho(\Lambda h_3,l)h_1,h_2,h_5]_H+\Lambda[h_1,\rho(\Lambda h_3,l)h_2,h_5]_H\\
=&0,
\end{align*}
which imply that Eqs. \eqref{4.1} and \eqref{4.2} hold. Similarly, we can prove that Eqs.\eqref{4.3},  \eqref{4.4} and \eqref{4.5} are true. The proof is finished.
\end{proof}

\begin{prop}
 Let $H\stackrel{\Lambda_1}{\longrightarrow}L$ and  $H\stackrel{\Lambda_2}{\longrightarrow}L$  be two nonabelian embedding tensor 3-Lie algebras and  $(f_L,f_{H})$   a homomorphism from $H\stackrel{\Lambda_1}{\longrightarrow}L$ to $H\stackrel{\Lambda_2}{\longrightarrow}L$.
 Then the induced representation $(L;\mathfrak{l}_{\Lambda_1},\mathfrak{m}_{\Lambda_1},\mathfrak{r}_{\Lambda_1})$ of the descendent 3-Leibniz algebra $(H,[-,-,-]_{\Lambda_1})$ and the induced representation $(L;\mathfrak{l}_{\Lambda_2},\mathfrak{m}_{\Lambda_2},\mathfrak{r}_{\Lambda_2})$ of the descendent 3-Leibniz algebra $(H,[-,-,-]_{\Lambda_2})$ satisfying the following equations:
  \begin{align}
f_L(\mathfrak{l}_{\Lambda_1}(h_1,h_2,l))=&\mathfrak{l}_{\Lambda_2}(f_H(h_1),f_H(h_2),f_L(l)),\label{4.6}\\
f_L(\mathfrak{m}_{\Lambda_1}(h_1,l,h_2))=&\mathfrak{m}_{\Lambda_2}(f_H(h_1),f_L(l),f_H(h_2)),\label{4.7}\\
f_L(\mathfrak{r}_{\Lambda_1}(l, h_1,h_2))=&\mathfrak{r}_{\Lambda_2}(f_L(l),f_H(h_1),f_H(h_2)), \ \forall h_1,h_2\in H,l\in L. \label{4.8}
\end{align}
In other words,  the following diagrams commute:
 $$\aligned
\xymatrix{
L \ar[d]^{\mathfrak{l}_{\Lambda_1}(h_1,h_2,-)}\ar[rr]^{f_L} & & L \ar[d]^{\mathfrak{l}_{\Lambda_2}(f_H(h_1),f_H(h_2),-)}  \\
L \ar[rr]^{f_L} && L,}
\endaligned
\aligned
\xymatrix{
L \ar[d]^{\mathfrak{m}_{\Lambda_1}(h_1,-,h_2)}\ar[rr]^{f_L} & & L \ar[d]^{\mathfrak{m}_{\Lambda_2}(f_H(h_1),-,f_H(h_2))}  \\
L \ar[rr]^{f_L} && L,}
\endaligned$$
$$
\aligned
\xymatrix{
L \ar[d]^{\mathfrak{r}_{\Lambda_1}(-,h_1,h_2)}\ar[rr]^{f_L} & & L \ar[d]^{\mathfrak{r}_{\Lambda_2}(-,f_H(h_1),f_H(h_2))}  \\
L \ar[rr]^{f_L} && L.}
\endaligned
$$
\end{prop}

\begin{proof} For all  $h_1,h_2\in H,l\in L,$ by Eqs. \eqref{2.7} and \eqref{2.8}
 we have
 \begin{align*}
f_L(\mathfrak{l}_{\Lambda_1}(h_1,h_2,l))=&f_L([\Lambda_1 h_1,\Lambda_1 h_2,l]_L)=[f_L(\Lambda_1 h_1),f_L(\Lambda_1 h_2),f_L(l)]_L\\
=&[\Lambda_2 f_H(h_1),\Lambda_2 f_H(h_2),f_L(l)]_L\\
=&\mathfrak{l}_{\Lambda_2}(f_H(h_1),f_H(h_2),f_L(l)),\\
f_L(\mathfrak{m}_{\Lambda_1}(h_1,l,h_2))=&f_L([\Lambda_1 h_1,l,\Lambda_1 h_2]_L-\Lambda_1\rho(\Lambda_1 h_1,l)h_2)\\
=&[f_L(\Lambda_1 h_1),f_L(l),f_L(\Lambda_1 h_2)]_L-f_L(\Lambda_1\rho(\Lambda_1 h_1,l)h_2)\\
=&[\Lambda_2 f_H(h_1),f_L(l),\Lambda_2 f_H(h_2)]_L-\Lambda_2f_H(\rho(\Lambda_1 h_1,l)h_2)\\
=&[\Lambda_2 f_H(h_1),f_L(l),\Lambda_2 f_H(h_2)]_L-\Lambda_2\rho(\Lambda_2 f_H(h_1),f_L(l))f_H(h_2)\\
=&\mathfrak{m}_{\Lambda_2}(f_H(h_1),f_L(l),f_H(h_2)).
\end{align*}
 And the other equation is similar to provable.
\end{proof}

For $n\geq 1$,  let  $\delta_\Lambda: \mathcal{C}^n_{\mathrm{3Leib}}(H,L)\rightarrow \mathcal{C}^{n+1}_{\mathrm{3Leib}}(H,L)$ be  the coboundary operator of the 3-Leibniz algebra $(H,[-,-,-]_\Lambda)$
 with coefficients in the representation $(L;\mathfrak{l}_\Lambda,\mathfrak{m}_\Lambda,\mathfrak{r}_\Lambda)$. More precisely, for all  $\phi\in \mathcal{C}^n_{\mathrm{3Leib}}(H,L), \mathfrak{H}_i=u_i\wedge v_i\in \wedge^2 H, 1\leq i\leq n$ and $w\in H$, we have
\begin{small}
\begin{align*}
&(\delta_\Lambda \phi)(\mathfrak{H}_1,\mathfrak{H}_2, \cdots, \mathfrak{H}_n,w)\\
=&\sum_{1\leq j<k\leq n}(-1)^j\phi(\mathfrak{H}_1,\cdots,\widehat{\mathfrak{H}_j},\cdots,\mathfrak{H}_{k-1},u_k \wedge[u_j,v_j,v_k]_{\Lambda}+[u_j,v_j,u_k]_{\Lambda}\wedge  v_k,\cdots,\mathfrak{H}_n,w)\\
&+\sum_{j=1}^n(-1)^j\phi(\mathfrak{H}_1,\cdots,\widehat{\mathfrak{H}_j},\cdots,\mathfrak{H}_{n},[u_j,v_j,w]_{\Lambda})\\
&+\sum_{j=1}^n(-1)^{j+1}\mathfrak{l}_\Lambda(\mathfrak{H}_j,\phi(\mathfrak{H}_1,\cdots,\widehat{\mathfrak{H}_j},\cdots,\mathfrak{H}_{n},w))\\
&+(-1)^{n+1}(\mathfrak{m}_\Lambda(u_n, \phi(\mathfrak{H}_1,\cdots,\mathfrak{H}_{n-1},v_n),w)+\mathfrak{r}_\Lambda(\phi(\mathfrak{H}_1,\cdots,\mathfrak{H}_{n-1},u_n), v_n, w)).
\end{align*}
\end{small}
In particular, for $\phi\in \mathcal{C}^1_{\mathrm{3Leib}}(H,L):= \mathrm{Hom}(H,L)$ and $u_1,v_1,w\in H,$
we have
\begin{align*}
(\delta_\Lambda \phi)(u_1,v_1,w)=&-\phi([u_1,v_1,w]_\Lambda)+\mathfrak{l}_\Lambda(u_1,v_1,\phi(w))+\mathfrak{m}_\Lambda(u_1,\phi(v_1),w)+\mathfrak{r}_\Lambda(\phi(u_1),v_1,w)\\
=&-\phi([u_1,v_1,w]_\Lambda)+[\Lambda u_1,\Lambda v_1,\phi(w)]_L+[\Lambda u_1,\phi(v_1),\Lambda w]_L\\
&-\Lambda\rho(\Lambda u_1,\phi(v_1))w+[\phi(u_1),\Lambda v_1,\Lambda w]_L-\Lambda\rho(\phi(u_1),\Lambda v_1)w.
\end{align*}

For any  $(a_1,a_2)\in  \mathcal{C}^0_{\mathrm{3Leib}}(V,L):=\wedge^2 L$, we define $\delta_\Lambda: \mathcal{C}^0_{\mathrm{3Leib}}(V,L)\rightarrow \mathcal{C}^1_{\mathrm{3Leib}}(V,L), (a_1, a_2)\mapsto \delta_\Lambda(a_1,a_2)$ by
$$\delta_\Lambda(a_1, a_2)u=\Lambda\rho(a_1, a_2)u-[a_1, a_2,\Lambda u]_L, \forall  u\in H.$$

\begin{prop}
Let $H\stackrel{\Lambda}{\longrightarrow}L$    be a nonabelian embedding tensor 3-Lie algebra. Then $\delta_\Lambda(\delta_\Lambda(a,b))=0$.
\end{prop}

\begin{proof}  For any $u_1,v_1,w\in V,$ by Eqs. \eqref{2.1}-\eqref{2.6} and \eqref{2.9} we have
\begin{align*}
&\delta_\Lambda (\delta_\Lambda(a_1,a_2))(u_1,v_1,w)\\
=&-\delta_\Lambda(a_1,a_2)([u_1,v_1,w]_\Lambda)+[\Lambda u_1,\Lambda v_1,\delta_\Lambda(a_1,a_2)(w)]_L+[\Lambda u_1,\delta_\Lambda(a_1,a_2)(v_1),\Lambda w]_L\\
&-\Lambda\rho(\Lambda u_1,\delta_\Lambda(a_1,a_2)(v_1))w+[\delta_\Lambda(a_1,a_2)(u_1),\Lambda v_1,\Lambda w]_L-\Lambda\rho(\delta_\Lambda(a_1,a_2)(u_1),\Lambda v_1)w\\
=&-\Lambda\rho(a_1, a_2)[u_1,v_1,w]_\Lambda+[a_1, a_2, [\Lambda u_1,\Lambda v_1,\Lambda w]_L]_L+[\Lambda u_1,\Lambda v_1,\Lambda\rho(a_1, a_2)w]_L\\
&-[\Lambda u_1,\Lambda v_1,[a_1, a_2,\Lambda w]_L]_L+[\Lambda u_1,\Lambda\rho(a_1, a_2)v_1,\Lambda w]_L-[\Lambda u_1,[a_1, a_2,\Lambda v_1]_L,\Lambda w]_L\\
&-\Lambda\rho(\Lambda u_1,\Lambda\rho(a_1, a_2)v_1)w+\Lambda\rho(\Lambda u_1,[a_1, a_2,\Lambda v_1]_L)w+[\Lambda\rho(a_1, a_2)u_1,\Lambda v_1,\Lambda w]_L\\
&-[[a_1, a_2,\Lambda u_1]_L,\Lambda v_1,\Lambda w]_L-\Lambda\rho(\Lambda\rho(a_1, a_2)u_1,\Lambda v_1)w+\Lambda\rho([a_1, a_2,\Lambda u_1]_L,\Lambda v_1)w\\
=&-\Lambda\rho(a_1, a_2)\rho(\Lambda u_1, \Lambda v_1)w-\Lambda\rho(a_1, a_2)[u_1,v_1,w]_H+\Lambda\rho(\Lambda u_1,\Lambda v_1)\rho(a_1, a_2)w\\
&+\Lambda[u_1, v_1, \rho(a_1, a_2)w]_H+\Lambda\rho(\Lambda u_1,\Lambda\rho(a_1, a_2)v_1) w+\Lambda[u_1, \rho(a_1, a_2)v_1, w]_H\\
&-\Lambda\rho(\Lambda u_1,\Lambda\rho(a_1, a_2)v_1)w+\Lambda\rho(\Lambda u_1,[a_1, a_2,\Lambda v_1]_L)w+\Lambda(\Lambda\rho(a_1, a_2)u_1,\Lambda v_1)w\\
&+\Lambda[\rho(a_1, a_2)u_1, v_1,w]_H-\Lambda\rho(\Lambda\rho(a_1, a_2)u_1,\Lambda v_1)w+\Lambda\rho([a_1, a_2,\Lambda u_1]_L,\Lambda v_1)w\\
=&-\Lambda\rho(a_1, a_2)\rho(\Lambda u_1, \Lambda v_1)w+\Lambda\rho(\Lambda u_1,\Lambda v_1)\rho(a_1, a_2)w+\Lambda\rho(\Lambda u_1,\Lambda\rho(a_1, a_2)v_1) w\\
&-\Lambda\rho(\Lambda u_1,\Lambda\rho(a_1, a_2)v_1)w+\Lambda\rho(\Lambda u_1,[a_1, a_2,\Lambda v_1]_L)w+\Lambda(\Lambda\rho(a_1, a_2)u_1,\Lambda v_1)w\\
&-\Lambda\rho(\Lambda\rho(a_1, a_2)u_1,\Lambda v_1)w+\Lambda\rho([a_1, a_2,\Lambda u_1]_L,\Lambda v_1)w\\
=&0.
\end{align*}
Therefore, we deduce that $\delta_\Lambda (\delta_\Lambda(a_1,a_2))=0.$
\end{proof}

\medskip

Now we develop the cohomology theory of a  nonabelian embedding tensor $\Lambda$ on the 3-Lie algebra     $(L, [-,-,-]_L)$ with
respect to the coherent action $(H,[-,-,-]_H;\rho^{\dag})$.

For $n\geq 0$,  define the set of $n$-cochains of $\Lambda$ by
$\mathcal{C}^n_\Lambda(H,L):=\mathcal{C}^n_{\mathrm{3Leib}}(H,L).$  Then $(\oplus_{n=0}^{\infty}\mathcal{C}^n_\Lambda(H,L),\delta_\Lambda)$ is a cochain complex.

For $n\geq 1$,  we denote the set of $n$-cocycles by
$\mathcal{Z}^n_\Lambda(H,L)$, the set of $n$-coboundaries by $\mathcal{B}^n_\Lambda(H,L)$ and the $n$-th cohomology group of the nonabelian embedding tensor $\Lambda$ by $\mathcal{H}^n_\Lambda(H,L)=\mathcal{Z}^n_\Lambda(H,L)/\mathcal{B}^n_\Lambda(H,L)$.

\begin{prop}
 Let $H\stackrel{\Lambda_1}{\longrightarrow}L$ and  $H\stackrel{\Lambda_2}{\longrightarrow}L$  be two nonabelian embedding tensor 3-Lie algebras and let  $(f_L,f_{H})$ be  a homomorphism from $H\stackrel{\Lambda_1}{\longrightarrow}L$ to $H\stackrel{\Lambda_2}{\longrightarrow}L$ in which
$f_H$ is invertible.
 We define a map $\Psi:\mathcal{C}^{n}_{\Lambda_1}(H,L)\rightarrow \mathcal{C}^{n}_{\Lambda_2}(H,L)$ by
   \begin{align*}
\Psi (\phi)(\mathfrak{H}_1,\mathfrak{H}_2, \cdots, \mathfrak{H}_{n-1},w)=f_L\big(\phi(f_H^{-1}(u_1)\wedge f_H^{-1}(v_1),\cdots, f_H^{-1}(u_{n-1})\wedge f_H^{-1}(v_{n-1}),f_H^{-1}(w) )\big),
\end{align*}
 for any $\phi\in \mathcal{C}^{n}_{\mathrm{\Lambda_1}}(H,L), \mathfrak{H}_i=u_i\wedge v_i\in \wedge^2 H, 1\leq i\leq {n-1}$ and $w\in H$.
 Then $\Psi: (\mathcal{C}^{n+1}_{\mathrm{\Lambda_1}}(H,L),\delta_{\Lambda_1})\rightarrow (\mathcal{C}^{n+1}_{\mathrm{\Lambda_2}}(H,L),\delta_{\Lambda_2})$  is a cochain map.

That is,  the following diagram commutes:
 $$\aligned
\xymatrix{
\mathcal{C}^n_{\mathrm{\Lambda_1}}(H,L) \ar[d]^{\Psi}\ar[rr]^{\delta_{\Lambda_1}} & & \mathcal{C}^{n+1}_{\mathrm{\Lambda_1}}(H,L) \ar[d]^{\Psi}  \\
\mathcal{C}^n_{\mathrm{\Lambda_2}}(H,L) \ar[rr]^{\delta_{\Lambda_2}} && \mathcal{C}^{n+1}_{\mathrm{\Lambda_2}}(H,L).}
\endaligned
$$
Consequently, it induces a homomorphism $\Psi^*$ from
the cohomology group $\mathcal{H}^{n+1}_{\mathrm{\Lambda_1}}(H,L)$   to  $\mathcal{H}^{n+1}_{\mathrm{\Lambda_2}}(H,L)$.
\end{prop}

\begin{proof}  For any  $\phi\in \mathcal{C}^{n}_{\mathrm{\Lambda_1}}(H,L), \mathfrak{H}_i=u_i\wedge v_i\in \wedge^2 H, 1\leq i\leq {n}$ and $w\in H$, by Eqs. \eqref{4.6}-\eqref{4.8} and  Proposition \ref{pro:2.9}, we have
\begin{align*}
 &(\delta_{\Lambda_2}\Psi(\phi))(\mathfrak{H}_1,\mathfrak{H}_2, \cdots, \mathfrak{H}_n,w)\\
=&\sum_{1\leq j<k\leq n}(-1)^j\Psi(\phi)(\mathfrak{H}_1,\cdots,\widehat{\mathfrak{H}_j},\cdots,\mathfrak{H}_{k-1},u_k \wedge[u_j,v_j,v_k]_{\Lambda_2}+[u_j,v_j,u_k]_{\Lambda_2}\wedge  v_k,\cdots,\mathfrak{H}_n,w)\\
&+\sum_{j=1}^n(-1)^j\Psi(\phi)(\mathfrak{H}_1,\cdots,\widehat{\mathfrak{H}_j},\cdots,\mathfrak{H}_{n},[u_j,v_j,w]_{\Lambda_2})\\
&+\sum_{j=1}^n(-1)^{j+1}\mathfrak{l}_{\Lambda_2}(\mathfrak{H}_j,\Psi(\phi)(\mathfrak{H}_1,\cdots,\widehat{\mathfrak{H}_j},\cdots,\mathfrak{H}_{n},w))\\
&+(-1)^{n+1}\mathfrak{m}_{\Lambda_2}(u_n, \Psi(\phi)(\mathfrak{H}_1,\cdots,\mathfrak{H}_{n-1},v_n),w)\\
&+(-1)^{n+1}\mathfrak{r}_{\Lambda_2}(\Psi(\phi)(\mathfrak{H}_1,\cdots,\mathfrak{H}_{n-1},u_n), v_n, w) \\
=&\sum_{1\leq j<k\leq n}(-1)^jf_L(\phi(f_H^{-1}(u_1)\wedge f_H^{-1}(v_1),\cdots, \widehat{\mathfrak{H}_j},\cdots,f_H^{-1}(u_{k-1})\wedge f_H^{-1}(v_{k-1}),\\
&f_H^{-1}(u_k) \wedge f_H^{-1}([u_j,v_j,v_k]_{\Lambda_2})+f_H^{-1}([u_j,v_j,u_k]_{\Lambda_2}) \wedge f_H^{-1}(v_k),f_H^{-1}(u_n)\wedge f_H^{-1}(v_n), f_H^{-1}(w) ))\\
&+\sum_{j=1}^n(-1)^jf_L(\phi(f_H^{-1}(u_1)\wedge f_H^{-1}(v_1),\cdots, \widehat{\mathfrak{H}_j},\cdots,f_H^{-1}(u_n)\wedge f_H^{-1}(v_n),f_H^{-1}([u_j,v_j,w]_{\Lambda_2})))\\
&+\sum_{j=1}^n(-1)^{j+1}\mathfrak{l}_{\Lambda_2}(\mathfrak{H}_j,f_L(\phi(f_H^{-1}(u_1)\wedge f_H^{-1}(v_1),\cdots, \widehat{\mathfrak{H}_j},\cdots,f_H^{-1}(u_n)\wedge f_H^{-1}(v_n),f_H^{-1}(w))))\\
&+(-1)^{n+1}\mathfrak{m}_{\Lambda_2}(u_n, f_L(\phi(f_H^{-1}(u_1)\wedge f_H^{-1}(v_1),\cdots, f_H^{-1}(u_{n-1})\wedge f_H^{-1}(v_{n-1}),f_H^{-1}(v_n))),w)\\
&+(-1)^{n+1}\mathfrak{r}_{\Lambda_2}(f_L(\phi(f_H^{-1}(u_1)\wedge f_H^{-1}(v_1),\cdots,f_H^{-1}(u_{n-1})\wedge f_H^{-1}(v_{n-1}),f_H^{-1}(u_n))), v_n, w)
\end{align*}
\begin{align*}
=&f_L\big(\sum_{1\leq j<k\leq n}(-1)^j\phi(f_H^{-1}(u_1)\wedge f_H^{-1}(v_1),\cdots, \widehat{\mathfrak{H}_j},\cdots,f_H^{-1}(u_{k-1})\wedge f_H^{-1}(v_{k-1}),\\
&f_H^{-1}(u_k) \wedge [f_H^{-1}(u_j),f_H^{-1}(v_j),f_H^{-1}(v_k)]_{\Lambda_1}+[f_H^{-1}(u_j),f_H^{-1}(v_j),f_H^{-1}(u_k)]_{\Lambda_1} \wedge f_H^{-1}(v_k),\\
&f_H^{-1}(u_n)\wedge f_H^{-1}(v_n), f_H^{-1}(w))+\sum_{j=1}^n(-1)^j\phi(f_H^{-1}(u_1)\wedge f_H^{-1}(v_1),\cdots, \widehat{\mathfrak{H}_j},\cdots,f_H^{-1}(u_n)\wedge f_H^{-1}(v_n),\\
&[f_H^{-1}(u_j),f_H^{-1}(v_j),f_H^{-1}(w)]_{\Lambda_1})+\sum_{j=1}^n(-1)^{j+1}\mathfrak{l}_{\Lambda_1}(f_H^{-1}(u_j), f_H^{-1}(v_j),\phi(f_H^{-1}(u_1)\wedge f_H^{-1}(v_1),\cdots, \\
&\widehat{\mathfrak{H}_j},\cdots,f_H^{-1}(u_n)\wedge f_H^{-1}(v_n),f_H^{-1}(w)))+(-1)^{n+1}\mathfrak{m}_{\Lambda_1}(f_H^{-1}(u_n), \phi(f_H^{-1}(u_1), f_H^{-1}(v_1),\cdots, \\
&f_H^{-1}(u_{n-1})\wedge f_H^{-1}(v_{n-1}),f_H^{-1}(v_n)),f_H^{-1}(w))+(-1)^{n+1}\mathfrak{r}_{\Lambda_1}(\phi(f_H^{-1}(u_1)\wedge f_H^{-1}(v_1),\cdots,\\
&f_H^{-1}(u_{n-1})\wedge f_H^{-1}(v_{n-1}),f_H^{-1}(u_n)),  f_H^{-1}(v_n), f_H^{-1}(w))\big)\\
=&f_L(\delta_{\Lambda_1}\phi)(f_H^{-1}(u_1)\wedge f_H^{-1}(v_1),\cdots,f_H^{-1}(u_{n})\wedge f_H^{-1}(v_{n}),f_H^{-1}(w))\\
=&\Psi(\delta_{\Lambda_1}\phi)(\mathfrak{H}_1,\mathfrak{H}_2, \cdots, \mathfrak{H}_n,w).
\end{align*}
Hence, $\Psi$  is a cochain map, and induces   a cohomology group homomorphism $\Psi^*: \mathcal{H}^{n+1}_{\mathrm{\Lambda_1}}(H,L)$ $\rightarrow \mathcal{H}^{n+1}_{\mathrm{\Lambda_2}}(H,L)$.
\end{proof}

At the end of this  section, we use the established cohomology theory to characterize infinitesimal deformations of nonabelian embedding tensors on 3-Lie algebras.

\begin{defn}
 Let $\Lambda: H\rightarrow L$   be a nonabelian embedding tensor  on a 3-Lie algebra $(L, [-,-,-]_L)$ with
respect to a coherent action $(H,[-,-,-]_H;\rho^\dag)$.  An infinitesimal deformation of $\Lambda$ is a  nonabelian embedding tensor of the form $\Lambda_t=\Lambda+t\Lambda_1$, where $t$ is a parameter
with $t^2=0.$
\end{defn}

Let $\Lambda+t\Lambda_1$ be an infinitesimal deformation of $\Lambda$, then we
have
\begin{align*}
[\Lambda_tu_1,\Lambda_tu_2, \Lambda_tu_3 ]_L=& \Lambda_t\rho(\Lambda_tu_1,\Lambda_tu_2)u_3+\Lambda_t[u_1,u_2,u_3]_H,
\end{align*}
for any $u_1,u_2,u_3\in H.$  This is equivalent to the following conditions
\begin{align}
&[\Lambda_1u_1,\Lambda u_2,\Lambda u_3]_L+[\Lambda u_1,\Lambda_1 u_2,\Lambda u_3]_L+[\Lambda u_1,\Lambda u_2,\Lambda_1 u_3]_L \nonumber\\
&=\Lambda_1\rho(\Lambda u_1,\Lambda u_2)u_3+\Lambda\rho(\Lambda_1 u_1,\Lambda u_2)u_3+\Lambda\rho(\Lambda u_1,\Lambda_1 u_2)u_3+\Lambda_1 [u_1,u_2,u_3]_H, \label{4.9}\\
&[\Lambda_1u_1,\Lambda_1 u_2,\Lambda u_3]_L+[\Lambda_1 u_1,\Lambda u_2,\Lambda_1 u_3]_L+[\Lambda u_1,\Lambda_1 u_2,\Lambda_1 u_3]_L \nonumber\\
&=\Lambda_1\rho(\Lambda_1 u_1,\Lambda u_2)u_3+\Lambda_1\rho(\Lambda u_1,\Lambda_1 u_2)u_3+\Lambda\rho(\Lambda_1 u_1,\Lambda_1 u_2)u_3, \label{4.10}\\
&[\Lambda_1u_1,\Lambda_1 u_2,\Lambda_1 u_3]_L=\Lambda_1\rho(\Lambda_1 u_1,\Lambda_1 u_2)u_3.\label{4.11}
\end{align}
 From Eq. \eqref{4.11}   it follows that the
map $\Lambda_1$ is  an   embedding tensor  on  the 3-Lie algebra $(L, [-,-,-])$ with
respect to the representation $(H;\rho)$ (see \cite{Hu}).
It follows from Eq. \eqref{4.9} that $\Lambda_1\in \mathcal{C}_\Lambda^1(H,L)$
is a 1-cocycle in the cohomology complex of $\Lambda$.
Thus the cohomology class of $\Lambda_1$ defines an element in $\mathcal{H}_\Lambda^1(H,L)$.

Let $\Lambda_t=\Lambda+t\Lambda_1$ and $\Lambda'_t=\Lambda+t\Lambda'_1$
 be two infinitesimal deformations of $\Lambda$. They are said to be equivalent if there
exist  $a_1\wedge a_2\in \wedge^2 L$ such that the pair  $(id_L+tad(a_1,a_2),id_H+t\rho(a_1,a_2))$
is a homomorphism   from $H\stackrel{\Lambda_t}{\longrightarrow}L$ to $H\stackrel{\Lambda'_t}{\longrightarrow}L$.
That is, the following
conditions must hold:

(1) the   maps $ id_L+tad(a_1,a_2):L\rightarrow L $ and  $ id_H+t\rho(a_1,a_2): H\rightarrow H $  are two 3-Lie algebra homomorphisms,

(2)  the pair   $(id_L+tad(a_1,a_2),id_H+t\rho(a_1,a_2))$ satisfies:
\begin{align}
&(id_H+t\rho(a_1,a_2))(\rho(a,b)u)\nonumber\\
&=\rho((id_L+tad(a_1,a_2))(a),(id_L+tad(a_1,a_2))(b))(id_H+t\rho(a_1,a_2))(u),\label{4.12} \\
&(\Lambda+t\Lambda'_1)(id_H+t\rho(a_1,a_2))(u)=(id_L+tad(a_1,a_2))(\Lambda+t\Lambda_1)(u),\label{4.13}
\end{align}
for all $a,b\in L, u\in H. $ It is easy to see that the condition \eqref{4.13}  gives rise to
 $$\Lambda_1 u-\Lambda'_1u=\Lambda\rho(a_1,a_2)u-[a_1,a_2,\Lambda u]=\delta_\Lambda(a_1,a_2)u\in \mathcal{C}_\Lambda^1(H,L).$$
This shows that $\Lambda_1$ and $\Lambda'_1$ are cohomologous. Thus their cohomology classes are the same   in $\mathcal{H}_\Lambda^1(H,L)$.

Conversely, any 1-cocycle $\Lambda_1$ gives rise to the infinitesimal deformation $\Lambda+t\Lambda_1$. Furthermore, we have the following result.

\begin{theorem}
 Let $\Lambda: H\rightarrow L$   be a nonabelian embedding tensor  on a 3-Lie algebra $(L, [-,-,-]_L)$ with
respect to a coherent action $(H,[-,-,-]_H;\rho^{\dag})$.    Then there is a bijection between the set of all
equivalence classes of infinitesimal deformations of $\Lambda$ and the first cohomology group $\mathcal{H}_\Lambda^1(H,L)$.
\end{theorem}

\section{Nonabelian embedding tensors on 3-Lie algebras induced by  Lie algebras}
\def\theequation{\arabic{section}.\arabic{equation}}
\setcounter{equation} {0}

Motivated by the construction of 3-Lie algebras from Lie algebras. In this section, we provide and investigate nonabelian embedding tensors on 3-Lie algebras induced by  Lie algebras.
Recall from \cite{Bai} that given a  Lie algebra and a trace map one can construct a
3-Lie algebra. Let $(L,[-,-]_L)$ be a  Lie algebra and $L^*$ the dual of $L$. $\varsigma\in L^*$ is called a trace map  if it satisfies
$\varsigma([l_1,l_2]_L)=0$, for any $l_1,l_2\in L$. We
define the ternary bracket $[-,-,-]_{L_\varsigma}$ by
$$[l_1,l_2,l_3]_{L_\varsigma}=\varsigma(l_1)[l_2,l_3]_L+\varsigma(l_2)[l_3,l_1]_L+\varsigma(l_3)[l_1,l_2]_L, \ \forall l_1,l_2, l_3\in L.$$
This 3-Lie algebra is denoted by $L_\varsigma$.

A coherent action of a  Lie algebra $(L,[-,-]_L)$ on a Lie algebra $(H,[-,-]_H)$ is a Lie algebra homomorphism $\rho:L\rightarrow \mathrm{Der}(H)$ satisfies $[\rho(l)h_1,h_2]_H=0, \ \forall l\in L, h_1,h_2\in H.$ See \cite{Tang} for more details. We denote a coherent action of $(L,[-,-]_L)$ by $(H,[-,-]_{H}; \rho^\dag)$.

\begin{prop}
let $(H,[-,-]_{H}; \rho^\dag)$ be  a coherent action of  a  Lie algebra  $(L,[-,-]_L)$ and $\varsigma_L,\varsigma_H$ be two   trace maps,   that is two
linear maps satisfying
$$\varsigma_L([l_1,l_2]_L)=0,\varsigma_H([h_1, h_2]_H)=0, \ \ \forall l_1,l_2\in L, h_1,h_2\in H.$$
In order to simplify, $\varsigma_L$ and $\varsigma_H$ are denoted by the same symbol $\varsigma$.
Then $(H,[-,-,-]_\varsigma; \rho_\varsigma^\dag)$ is a coherent action of  the 3-Lie algebra   $L_\varsigma$,
 where $\rho_\varsigma:\wedge^2 L\rightarrow \mathrm{End}(H)$ is defined by $$\rho_\varsigma(l_1,l_2)=\varsigma(l_1)\rho(l_2)-\varsigma(l_2)\rho(l_1),\ \ \forall l_1,l_2\in L.$$
\end{prop}
\begin{proof}
In the light of \cite{Zhao}(Proposition 5.6),  $(H; \rho_\varsigma)$ is a  representation of  a 3-Lie algebra   $L_\varsigma$.  We only need to check that $\rho_\varsigma$ satisfies Eqs.  \eqref{2.4} and  \eqref{2.5}.
For all $l_1,l_2\in L$ and $h_1,h_2,h_3\in H$, we have
\begin{align*}
&[\rho_\varsigma(l_1,l_2)h_1,h_2,h_3]_{H_\varsigma}+[h_1,\rho_\varsigma(l_1,l_2)h_2,h_3]_{H_\varsigma}+[h_1,h_2,\rho_\varsigma(l_1,l_2)h_3]_{H_\varsigma}\\
=&[(\varsigma(l_1)\rho(l_2)-\varsigma(l_2)\rho(l_1))h_1,h_2,h_3]_{H_\varsigma}+[h_1,(\varsigma(l_1)\rho(l_2)-\varsigma(l_2)\rho(l_1))h_2,h_3]_{H_\varsigma}\\
&+[h_1,h_2,(\varsigma(l_1)\rho(l_2)-\varsigma(l_2)\rho(l_1))h_3]_{H_\varsigma}\\
=&\varsigma(l_1)\varsigma(\rho(l_2)h_1)[h_2,h_3]_H+\varsigma(l_1)\varsigma(h_2)[h_3,\rho(l_2)h_1]_H+\varsigma(l_1)\varsigma(h_3)[\rho(l_2)h_1,h_2]_H\\
&-\varsigma(l_2)\varsigma(\rho(l_1)h_1)[h_2,h_3]_H-\varsigma(l_2)\varsigma(h_2)[h_3,\rho(l_1)h_1]_H-\varsigma(l_2)\varsigma(h_3)[\rho(l_1)h_1,h_2]_H\\
&+\varsigma(l_1)\varsigma(h_1)[\rho(l_2)h_2,h_3]_H+\varsigma(l_1)\varsigma(\rho(l_2)h_2)[h_3,h_1]_H+\varsigma(l_1)\varsigma(h_3)[h_1,\rho(l_2)h_2]_H\\
&-\varsigma(l_2)\varsigma(h_1)[\rho(l_1)h_2,h_3]_H-\varsigma(l_2)\varsigma(\rho(l_1)h_2)[h_3,h_1]_H-\varsigma(l_2)\varsigma(h_3)[h_1,\rho(l_1)h_2]_H\\
&+\varsigma(l_1)\varsigma(h_1)[h_2,\rho(l_2)h_3]_H+\varsigma(l_1)\varsigma(h_2)[\rho(l_2)h_3,h_1]_H+\varsigma(l_1)\varsigma(\rho(l_2)h_3)[h_1,h_2]_H\\
&-\varsigma(l_2)\varsigma(h_1)[h_2,\rho(l_1)h_3]_H-\varsigma(l_2)\varsigma(h_2)[\rho(l_1)h_3,h_1]_H-\varsigma(l_2)\varsigma(\rho(l_1)h_3)[h_1,h_2]_H\\
=&\varsigma(l_1)\varsigma(\rho(l_2)h_1)[h_2,h_3]_H-\varsigma(l_2)\varsigma(\rho(l_1)h_1)[h_2,h_3]_H+\varsigma(l_1)\varsigma(\rho(l_2)h_2)[h_3,h_1]_H\\
&-\varsigma(l_2)\varsigma(\rho(l_1)h_2)[h_3,h_1]_H+\varsigma(l_1)\varsigma(\rho(l_2)h_3)[h_1,h_2]_H-\varsigma(l_2)\varsigma(\rho(l_1)h_3)[h_1,h_2]_H\\
=&(\varsigma(l_1)\rho(l_2)-\varsigma(l_2)\rho(l_1))(\varsigma(h_1)[h_2,h_3]_H+\varsigma(h_2)[h_3,h_1]_H+\varsigma(h_3)[h_1,h_2]_H)\\
=&\rho_\varsigma(l_1,l_2)[h_1,h_2,h_3]_{H_\varsigma},\\
&[\rho_\varsigma(l_1,l_2)h_1,h_2,h_3]_{H_\varsigma}\\
=&[\varsigma(l_1)\rho(l_2)h_1-\varsigma(l_2)\rho(l_1)h_1,h_2,h_3]_{H_\varsigma} \\
=&\varsigma(l_1)[\rho(l_2)h_1,h_2,h_3]_{H_\varsigma}-\varsigma(l_2)[\rho(l_1)h_1,h_2,h_3]_{H_\varsigma} \\
=& \varsigma(l_1)\varsigma(\rho(l_2)h_1)[h_2,h_3]_H+\varsigma(l_1)\varsigma(h_2)[h_3,\rho(l_2)h_1]_H+\varsigma(l_1)\varsigma(h_3)[\rho(l_2)h_1,h_2]_H\\
&-\varsigma(l_2)\varsigma(\rho(l_1)h_1)[h_2,h_3]_H-\varsigma(l_2)\varsigma(h_2)[h_3,\rho(l_1)h_1]_H-\varsigma(l_2)\varsigma(h_3)[\rho(l_1)h_1,h_2]_H\\
=&0.
\end{align*}
The proof is finished.
\end{proof}

A nonabelian embedding tensor on a Lie algebra $(L,[-,-]_L)$ with respect to a
coherent action $(H,[-,-]_{H}; \rho^\dag)$ is a linear map $\Lambda: H\rightarrow L$  such that
$[\Lambda h_1,\Lambda h_2]_L=\Lambda(\rho(\Lambda h_1)h_2+[h_1,h_2]_H), \forall h_1,h_2\in H$ (see \cite{Tang}).

\begin{theorem}
Let $\Lambda: H\rightarrow L$ be a nonabelian embedding tensor on a Lie algebra
$(L,[-,-]_L)$ with respect to a coherent action $(H,[-,-]_H;\rho^{\dag})$, let $\varsigma_L,\varsigma_H$ be two trace maps and  satisfying $\varsigma_L(\Lambda h)=\varsigma_H(h)$, for any $h\in H$.
In order to simplify, $\varsigma_L$ and $\varsigma_H$ are denoted by the same symbol $\varsigma$.
Then $\Lambda: H\rightarrow L$ is a nonabelian embedding tensor on the 3-Lie algebra
$(L,[-,-,-]_{L_\varsigma})$ with respect to the coherent action $(H,[-,-,-]_{H_\varsigma};\rho_{\varsigma}^{\dag})$.
\end{theorem}

\begin{proof}
For all   $h_1,h_2,h_3\in H$, we have
\begin{align*}
&[\Lambda h_1,\Lambda h_2, \Lambda h_3]_{L_\varsigma}\\
=&\varsigma(\Lambda h_1)[\Lambda h_2,\Lambda h_3]_L+\varsigma(\Lambda h_2)[\Lambda h_3,\Lambda h_1]_L+\varsigma(\Lambda h_3)[\Lambda h_1,\Lambda h_2]_L\\
=&\varsigma(\Lambda h_1)(\Lambda \rho(\Lambda h_2)h_3+\Lambda[h_2,h_3]_H)+\varsigma(\Lambda h_2)(\Lambda \rho(\Lambda h_3)h_1+\Lambda[h_3,h_1]_H)\\
&+\varsigma(\Lambda h_3)(\Lambda \rho(\Lambda h_1)h_2+\Lambda[h_1,h_2]_H)\\
=&\varsigma(\Lambda h_1)\Lambda \rho(\Lambda h_2)h_3+\varsigma(\Lambda h_2)\Lambda \rho(\Lambda h_3)h_1+\varsigma(\Lambda h_3)\Lambda \rho(\Lambda h_1)h_2+\varsigma(\Lambda h_2)\Lambda[h_3,h_1]_H\\
&+\varsigma(\Lambda h_1)\Lambda[h_2,h_3]_H+\varsigma(\Lambda h_3)\Lambda[h_1,h_2]_H.
\end{align*}
On the other hand, we have
\begin{align*}
&\Lambda(\rho_\varsigma(\Lambda h_1,\Lambda h_2)h_3+[h_1,h_2,h_3]_{H_\varsigma})\\
=&\Lambda (\varsigma(\Lambda h_1)\rho(\Lambda h_2)h_3- \varsigma(\Lambda h_2)\rho(\Lambda h_1)h_3+\varsigma( h_1)[h_2,h_3]_H+\varsigma( h_2)[h_3,h_1]_H+\varsigma( h_3)[h_1,h_2]_H).
\end{align*}
Thus, $[\Lambda h_1,\Lambda h_2, \Lambda h_3]_{L_\varsigma}=\Lambda(\rho_\varsigma(\Lambda h_1,\Lambda h_2)h_3+[h_1,h_2,h_3]_{H_\varsigma}),$ which implies that $\Lambda: H\rightarrow L$ is a nonabelian embedding tensor on
$(L,[-,-,-]_{L_\varsigma})$ with respect to  $(H,[-,-,-]_{H_\varsigma};\rho_{\varsigma}^{\dag})$.
\end{proof}

\begin{defn} (See \cite{Tang})
A Leibniz-Lie algebra $(H,[-,-]_H,\rhd)$ consists of a Lie algebra
$(H,[-,-]_H)$ and a binary product $\rhd:H\otimes H\rightarrow H$ such that
\begin{align*}
h_1\rhd (h_2\rhd h_3)&=(h_1\rhd  h_2)\rhd h_3+h_2\rhd (h_1\rhd h_3)+[h_1,h_2]_H\rhd h_3\\
h_1\rhd [h_2, h_3]_H&=[h_1\rhd h_2,h_3]_H=0,
\end{align*}
for all $h_1,h_2,h_3\in H$.
\end{defn}

\begin{theorem}
Let $(H,[-,-]_H,\rhd)$  be a Leibniz-Lie algebra  and $\varsigma$ be a trace map, that is a
linear map satisfying
$$\varsigma([h_1,h_2]_H)=0,\varsigma(h_1\rhd h_2)=0, \ \ \forall h_1,h_2\in H.$$
 Define two  ternary brackets by
\begin{align*}
 [h_1,h_2,h_3]_{H_\varsigma}&=\varsigma(h_1)[h_2,h_3]_H+\varsigma(h_2)[h_3,h_1]_H+\varsigma(h_3)[h_1,h_2]_H,\\
 \{h_1,h_2,h_3\}_{H_\varsigma}&=\varsigma(h_1) h_2\rhd h_3 -\varsigma(h_2) h_1\rhd h_3,  \forall h_1,h_2, h_3\in H.
\end{align*}
Then $(H,[-,-,-]_{H_\varsigma},\{-,-,-\}_{H_\varsigma})$ is a 3-Leibniz-Lie algebra.
\end{theorem}

\begin{proof}
For any $h_1,h_2,h_3,h_4,h_5\in H$, we have
\begin{align*}
&\{\{h_1,h_2,h_3\}_{H_\varsigma},h_4,h_5\}_{H_\varsigma}+\{h_3,\{h_1,h_2,h_4\}_{H_\varsigma}, h_5\}_{H_\varsigma} +\{h_3,h_4,\{h_1,h_2,h_5\}_{H_\varsigma}\}_{H_\varsigma} \\
&+\{[h_1,h_2,h_3]_{H_\varsigma},h_4,h_5\}_{H_\varsigma} +\{h_3,[h_1,h_2,h_4]_{H_\varsigma}, h_5\}_{H_\varsigma}-\{h_1,h_2,\{h_3,h_4,h_5\}_{H_\varsigma}\}_{H_\varsigma} \\
=&\varsigma(h_1)\varsigma( h_2\rhd h_3) h_4\rhd h_5 -\varsigma(h_4) \varsigma(h_1) (h_2\rhd h_3)\rhd h_5-\varsigma(h_2) \varsigma(h_1\rhd h_3) h_4\rhd h_5 \\
&+\varsigma(h_4) \varsigma(h_2) (h_1\rhd h_3)\rhd h_5+\varsigma(h_3) \varsigma(h_1) (h_2\rhd h_4)\rhd h_5 -\varsigma(h_1)\varsigma( h_2\rhd h_4) h_3\rhd h_5\\
&-\varsigma(h_3) \varsigma(h_2) (h_1\rhd h_4)\rhd h_5 +\varsigma(h_2)\varsigma( h_1\rhd h_4) h_3\rhd h_5+\varsigma(h_1)\varsigma(h_3) h_4\rhd( h_2\rhd h_5)\\
&-\varsigma(h_1)\varsigma(h_4) h_3\rhd (h_2\rhd h_5)-\varsigma(h_2)\varsigma(h_3) h_4\rhd( h_1\rhd h_5) +\varsigma(h_2)\varsigma(h_4) h_3\rhd (h_1\rhd h_5)\\
&+\varsigma(h_1)\varsigma([h_2,h_3]_H) h_4\rhd h_5 -\varsigma(h_4) \varsigma(h_1)[h_2,h_3]_H\rhd h_5+\varsigma(h_2)\varsigma([h_3,h_1]_H) h_4\rhd h_5 \\
&-\varsigma(h_4)\varsigma(h_2)[h_3,h_1]_H\rhd h_5+\varsigma(h_3)\varsigma([h_1,h_2]_H) h_4\rhd h_5 -\varsigma(h_4) \varsigma(h_3)[h_1,h_2]_H\rhd h_5\\
&+\varsigma(h_3) \varsigma(h_1)[h_2,h_4]_H\rhd h_5 -\varsigma(h_1)\varsigma([h_2,h_4]_H) h_3\rhd h_5+\varsigma(h_3) \varsigma(h_2)[h_4,h_1]_H\rhd h_5 \\
&-\varsigma(h_2)\varsigma([h_4,h_1]_H )h_3\rhd h_5+\varsigma(h_3) \varsigma(h_4)[h_1,h_2]_H\rhd h_5 -\varsigma(h_4)\varsigma([h_1,h_2]_H) h_3\rhd h_5\\
&-\varsigma(h_1)\varsigma(h_3) h_2\rhd (h_4\rhd h_5) +\varsigma(h_2)\varsigma(h_3) h_1\rhd (h_4\rhd h_5)+\varsigma(h_1)\varsigma(h_4) h_2\rhd (h_3\rhd h_5) \\
&-\varsigma(h_2)\varsigma(h_4) h_1\rhd (h_3\rhd h_5)\\
=&0,\\
&\{h_1,h_2,[h_3,h_4,h_5]_{H_\varsigma}\}_{H_\varsigma}\\
=&\varsigma(h_1) \varsigma(h_3)h_2\rhd[h_4,h_5]_H -\varsigma(h_2) \varsigma(h_3) h_1\rhd[h_4,h_5]_H+\varsigma(h_1)\varsigma(h_4) h_2\rhd[h_5,h_3]_H \\
&-\varsigma(h_2) \varsigma(h_4)h_1\rhd[h_5,h_3]_H+\varsigma(h_1)\varsigma(h_5) h_2\rhd[h_3,h_4]_H -\varsigma(h_2) \varsigma(h_5)h_1\rhd[h_3,h_4]_H\\
=&0,\\
&[\{h_1,h_2,h_3\}_{H_\varsigma},h_4,h_5]_{H_\varsigma}\\
=&\varsigma(h_1)\varsigma( h_2\rhd h_3)[h_4,h_5]_H+\varsigma(h_4)\varsigma(h_1)[h_5, h_2\rhd h_3]_H+\varsigma(h_5)\varsigma(h_1)[ h_2\rhd h_3,h_4]_H\\
&-\varsigma(h_2)\varsigma( h_1\rhd h_3)[h_4,h_5]_H-\varsigma(h_4)\varsigma(h_2)[h_5, h_1\rhd h_3]_H-\varsigma(h_5)\varsigma(h_2)[ h_1\rhd h_3,h_4]_H\\
=&0.
\end{align*}
The proof is finished.
\end{proof}

{{\bf Acknowledgments.}  The paper is  supported by the  Foundation of Science and Technology of Guizhou Province(Grant Nos. [2018]1020,  ZK[2022]031, ZK[2023]025),  the National Natural Science Foundation of China (Grant No. 12161013).


\begin{thebibliography}{99}


\bibitem{Arfa}  A. Arfa, N. Ben Fraj and A. Makhlouf, Cohomology and deformations of $n$-Lie algebra morphisms, {\it J. Geom. Phys.}, {\bf 132}, (2018), 64-74.



\bibitem{Bagger}  J. Bagger,  N. Lambert, Gauge symmetry and supersymmetry of multiple M2-branes, {\it Phys. Rev. D}, {\bf 77} (2008), no. 6, 065008.



\bibitem{Bai} R. Bai, Y  Wu, J Wang,  Realizatons of 3-Lie algebras. {\it J Math Phys.},  {\bf 51 }2010, Article ID 063305.

\bibitem{Bonezzi1}  R. Bonezzi, O. Hohm, Leibniz gauge theories and infinity structures, {\it Comm. Math. Phys.}, {\bf 377} (2020), 2027-2077.

\bibitem{Bonezzi2}  R. Bonezzi,  O. Hohm, Duality hierarchies and differential graded Lie algebras, {\it arXiv preprint}, 2019, arXiv:1910.10399.

\bibitem{Bergshoeff} E. A. Bergshoeff, M. de Roo,  O. Hohm,  Multiple M2-branes and the embedding tensor, {\it  Classical Quantum Gravity}, {\bf 25}(2008), 142001.


\bibitem{Casas} J. M. Casas,  J.  L. Loday,  T. Pirashvili, Leibniz $n$-algebras, {\it Forum Math.}, {\bf 14} (2002), 189-207.


\bibitem{deWit1}   B. de Wit, H. Samtleben,  M. Trigiante, On Lagrangians and gaugings of maximal supergravities, {\it Nuclear Phys. B},  {\bf 655} (2003),  1-2, 93-126.

\bibitem{deWit2}  B. de Wit, H. Samtleben, M. Trigiante, The maximal $D=5$ supergravities, {\it Nucl. Phys. B}, {\bf 716} (2005), 215-247.

\bibitem{deWit3}  B. de Wit,  H. Samtleben, Gauged maximal supergravities and hierarchies of nonabelian vector-tensor systems, {\it Fortschr. Phys.}, {\bf 53} (2005), 442-449.

\bibitem{deWit4}  B. de Wit, H. Nicolai, H. Samtleben, Gauged supergravities, tensor hierarchies, and M-theory, {\it J. High Energy Phys.}, {\bf 02} (2008), 044, 33 pp.



\bibitem{Das} A. Das, A. Makhlouf,  Embedding tensors on Hom-Lie algebras, {\it arXiv preprint}, 2023, arXiv:2304.04178.

\bibitem{Das1} A. Das,  Controlling structures, deformations and homotopy theory for averaging algebras,  {\it arXiv preprint}, 2023, arXiv:2303.17798.




\bibitem{Filippov}  V.T. Filippov, $n$-Lie algebras. {\it Sib. Mat. Zh.}, {\bf 26}, (1985), 126-140.


\bibitem{Gustavsson}  A. Gustavsson, Algebraic structures on parallel M2-branes, {\it Nuclear Phys. B}, {\bf 811} (2009), no. 1-2, 66-76.

\bibitem{Hu}  M. Hu, S. Hou, L. Song, Y. Zhou,  Deformations and cohomologies of embedding tensors on 3-Lie algebras, 2023, arXiv:2302.08725.

\bibitem{Ho} P. Ho, R. Hou,  Y. Matsuo, Lie 3-algebra and multiple M2-branes, {\it J. High Energy Phys.}  (2008) no. 6, 020,
30.

\bibitem{Hou} S. Hou, Y.  Sheng, Y. Zhou, 3-post-Lie algebras and relative Rota-Baxter operators of nonzero weight on 3-Lie algebras, {\it Journal of Algebra}, {\bf 615} 1,  (2023),   103-129.


\bibitem{Kasymov}  S.M. Kasymov, On a theory of $n$-Lie algebras (Russian). {\it Algebra i Logika}, {\bf 26}, (1987), 277-297.

\bibitem{Kotov}   A. Kotov, T. Strobl, The embedding tensor, Leibniz-Loday algebras and their higher gauge theories, {\it Comm. Math. Phys.},  {\bf 376} (2020), 235-258.



\bibitem{Liu}  J. Liu, Y. Sheng, Y. Zhou,  C. Bai, Nijenhuis operators on $n$-Lie algebras. {\it Commun. Theor. Phys.},  {\bf 65}, (2016), 659-670.

\bibitem{Liu1}  J. Liu, A. Makhlouf and Y. Sheng, A new approach to representations of 3-Lie algebras and abelian extensions. {\it Algebr. Represent. Theor.}  {\bf 20}, (2017), 1415-1431.

\bibitem{LiuW} W. Liu,  Z. Zhang, $T^*$-extension of a 3-Lie algebra, {\it Linear and Multilinear Algebra}, {\bf 60}5, (2012), 583-594


\bibitem{Mignel} J. Mignel, F. O. Farrill, Deformations of 3-algebras, {\it J. Math. Phys.}, {\bf 50},  (2009), 113514.



\bibitem{Nambu} Y. Nambu, Generalized Hamiltonian dynamics, {\it Phys. Rev. D}, {\bf 7} (1973), 2405-2412.

\bibitem{Nicolai}  H. Nicolai, H. Samtleben,  Maximal gauged supergravity in three dimensions, {\it Phys. Rev. Lett.},  {\bf 86} (2001),
1686-1689.



\bibitem{Sheng}  Y. Sheng, R. Tang, C. Zhu,   The controlling L$_{\infty}$-algebra, cohomology and homotopy of embedding tensors
and Lie-Leibniz triples, {\it Comm. Math. Phys.}, {\bf 386} (2021), 269-304.

\bibitem{Sheng1} Y. Sheng,  R. Tang,  Symplectic, product and complex structures on 3-Lie algebras,  {\it Journal of Algebra}, {\bf 208}, (2018), 256-300.

\bibitem{Takhtajan} L. Takhtajan, Higher order analog of Chevalley-Eilenberg complex and deformation theory of $n$-algebras, {\it St. Petersburg Math. J.}, {\bf 6}, (1995), 429-438.


\bibitem{Tang}  R. Tang, Y. Sheng,  Nonabelian embedding tensors,  {\it Lett. Math. Physics}, {\bf 113}(1)(2023),1-14.

\bibitem{Takhtajan}  L. A. Takhtajan, Higher order analog of Chevalley-Eilenberg complex and deformation theory of $n$-gebras, {\it St. Petersburg Math. J.}, {\bf 6} (1995), 429-438.

\bibitem{Teng} W. Teng, J. Jin, Y. Zhang, Cohomology of nonabelian embedding tensors on Hom-Lie algebras, {\it AIMS Mathematics}, {\bf 8}(9),2023, 21176-21190.

\bibitem{Teng1} W. Teng, J. Jin, Y. Zhang, Embedding tensors on 3-Hom-Lie algebras, {\it arXiv preprint}, 2023,	arXiv:2307.00961.

\bibitem{Xu} S. Xu, Cohomology, derivtions and abelian extensions of 3-Lie algebras. {\it J. Algebra Appl.}, {\bf 7}, (2019), 1950130, 26 pp.

\bibitem{Zhang}  T. Zhang, Deformations and Extensions of 3-Lie algebras, {\it arXiv preprint}, 2014,	arXiv:1401.4656.

\bibitem{Zhao} J. Zhao, J. Liu,  Y. Sheng, Cohomologies and relative Rota-Baxter-Nijenhuis structures of 3-LieRep pairs, {\it Linear and Multilinear Algebra},  {\bf 70}(21),(2021):1-25.


\end{thebibliography}
\end{document}